\documentclass[11pt]{article}
\input amssym.def
\usepackage{graphicx}

\def\sqw{\hbox{\rlap{\leavevmode\raise.3ex\hbox{$\sqcap$}}$%
\sqcup$}}
\def\sqb{\hbox{\vrule width5pt height5pt depth1.5pt%
}}

\def\ord{{\rm ord} }
  \def\bn{\hbox{\it 
I\hskip -2pt N}}      
\def\demo{\noindent{\bf Proof \ }} \newtheorem{theorem}{Theorem} 
\newtheorem{lemma}{Lemma} \newtheorem{proposition}{Proposition} 
\newtheorem{definition}{Definition} \newtheorem{example}{Example} 
\newtheorem{remark}{Remark}  
\newtheorem{corollary}{Corollary}
\begin{document}
\begin{center}
\uppercase{{\bf  The Nash problem on arcs for surface singularities }}
\end{center}
\advance\baselineskip-3pt
\vspace{2\baselineskip}
\begin{center}
{{\sc Marcel
Morales}\\
{\small Universit\'e de Grenoble I, Institut Fourier, 
UMR 5582, B.P.74,\\
38402 Saint-Martin D'H\`eres Cedex,\\
and IUFM de Lyon, 5 rue Anselme,\\ 69317 Lyon Cedex (FRANCE)}\\
 }
\vspace{\baselineskip}
\end{center}
\vskip.5truecm\noindent
{\small \sc Abstract.}{\footnote { first version  july 2005, revised version december 2005}} 
{\small Let $(X,O)$ be a germ of a normal surface singularity,  $\pi : \tilde X\longrightarrow X$ be the minimal 
resolution of singularities and let $A=(a_{ i,j} )$ be the $n\times 
n$ symmetrical intersection matrix of the exceptional  set of $\tilde X$. In an old 
preprint Nash proves that the set of arcs on a surface singularity is 
a scheme ${ \cal H }$, and defines a map ${ \cal N }$ from the set of irreducible 
components of ${ \cal H }$ to the set   of  
exceptional components of the minimal resolution of singularities of 
$(X,O)$.  He proved that this map is injective and ask if it  is 
surjective.
In this paper we consider the canonical decomposition 
 ${ \cal H }= \cup_{i=1}^n \bar{\cal N}_{i}$ :

\begin{itemize}

\item For any couple $(E_i,E_j)$ of distinct exceptional components, we define Numerical Nash
 condition $(NN_{(i,j)})$. We have that  $(NN_{(i,j)})$implies $ \bar{ 
\cal N}_{i}\not\subset \bar{ \cal N}_{j} $. In this paper   we  
 prove that $(NN_{(i,j)})$ is always true for at least the half 
of couples $(i,j)$. 
\item The condition $(NN_{(i,j)})$ is true for all couples $(i,j)$ with $i\not= j$, 
characterizes a certain class of negative definite  matrices, 
that we call Nash matrices. 
If $A$ is a Nash matrix then the Nash map ${ \cal N }$ is bijective. In particular
 our results 
depends only on $A$ and not on the topological type of the exceptional set. 
\item We recover and improve considerably almost all results known on this topic and our proofs 
are new  and  elementary. 
\item We give infinitely many other classes of singularities where Nash Conjecture is true.
\end{itemize}

The proofs are based on my old work \cite{M} and in Plenat \cite{P}.

\vskip 1cm 
\section{Introduction}
Let $(X,O)$ be a germ of a normal surface singularity. In an old 
preprint, published recently by Duke \cite{N}, Nash proved that the set of arcs on a surface singularity is 
a scheme $\cal{H}$, and defined a map ${ \cal N }$ from the set of irreducible 
components of $\cal{H}$ to the set   of  
exceptional components of the minimal resolution of singularities of 
$(X,O)$.  He proved that this map is injective and  ask if it  is 
surjective. 

Among the principal contributions to this subject we can cite  Monique 
Lejeune-Jalabert \cite{L-J}, Ana Reguera  \cite{R}, S. Ishii and J. Kollar \cite{I-K}, G. Gonzalez-Sprinberg and Monique Lejeune-Jalabert\cite{G-L}, Camille Plenat  \cite{P} and C. Plenat and P. Popescu-Pampu  \cite{P-PP}.
The study of arcs spaces was further developed by Kontsevich, Denef and Loeser \cite{D-L} in the theory of motivic integration.

Let $\pi : \tilde X\longrightarrow X$ be the minimal 
resolution of singularities, and $E_{1},\ldots,E_{n}$ be the 
components of the exceptional divisor, Ana Reguera  \cite{R} associates to every 
$E_{i}$ the family of arcs ${ \cal N}_{i}$ such that the proper 
transform cuts properly $E_{i}$, the spaces $\bar{ \cal N}_{i}$ are irreducible 
and give a decomposition of the space of arcs ${ \cal H }= \cup \bar{ 
\cal N}_{i}$.  In order to give an affirmative answer to the Nash 
problem it is sufficient to prove that for any $i\not=j $ then $ \bar{ 
\cal N}_{i}\not\subset \bar{ \cal N}_{j} $.

Recently Camille Plenat \cite{P}, Proposition 2.2 gives the following criterion to 
separate two Nash components:

\begin{proposition} Let  $\pi 
: \tilde X\longrightarrow X$ be the minimal resolution of 
singularities and $E_{1},\ldots,E_{n}$ be the components of the 
exceptional divisor, if there exist some $f\in { \cal O }_{X,O}$ such 
that $\ord _{E_{i} }(f)< \ord _{E_{j} }(f)$ then $ \bar{ \cal 
N}_{i}\not\subset \bar{ \cal N}_{j} $.
\end{proposition}

The following Theorem follows from  my work \cite{M} Theorem 1.1 and Lemma 2.2.
 Remark that in \cite{P-PP} C. Plenat and P. Popescu-Pampu have recently 
rediscover a similar condition.
\begin{theorem} Let $(X,O)$ be a germ of normal surface singularity,
 $\pi:\tilde X \longrightarrow X$  be the minimal
resolution of singularities and $E_{1},\ldots,E_{n}$ be the components 
of the exceptional divisor. Let $K_{\tilde X}$ the canonical divisor 
on ${\tilde X}$.  Let $E$ be an exceptional effective divisor and $Q= 
\pi_{*} { \cal O}_{\tilde X}¥ (-E) $,
\begin{enumerate}
\item If $-E\cdot E_{i} \geq 2K\cdot E_{i}$ for all $i=1,\ldots,n$ then 
$Q { \cal O}_{\tilde X}¥ ={ \cal O}_{\tilde X}¥ (-E)$

\item For any general  linear combination $f$ of a set of generators of 
$Q$ we have $div(f\circ \pi)= \tilde H+ E$, where $\tilde H$ is the proper 
transform of the cycle defined by $f$.
\end{enumerate}¥
\end{theorem}
\begin{remark} \begin{itemize}
\item  For any irreducible component $E_{i}$ of the exceptional divisor, we consider the adjunction
 formula for (eventually singular) curves 
$$p(E_i)=\frac{E_i\cdot (E_i+K_{\tilde X})}{2}+1$$ where $p(E_i)$ is the genus of $E_i$. Recall that $p(E_i)\geq 0$ and 
  $p(E_i)= 0$ if and only if $E_i$ is a curve of genus zero and self intersection equal to $-1$, 
which is impossible by Castelnuovo theorem  since we are assuming that $\pi:\tilde X \longrightarrow X$ is the minimal
resolution of singularities of $X$. As a consequence $K_{\tilde X}\cdot E_i=2(p(E_i)-1)-E_i^2\geq 0$ for any $i=1,...,n$.
\item  Since the graph of the resolution is connected we have that for any $1\leq i<k\leq n$ the intersection number 
$E_i\cdot E_k\geq 0$ and for each index $k$ there are at least one index $i$ such that  $E_i\cdot E_k> 0$.  
\item It follows from the previous item that if $\displaystyle E=\sum_{k=1}^{n} n_kE_{k}
 , n_k\in \bn $ is an exceptional divisor such that 
$E\cdot E_{k} \leq -2K_{\tilde X}\cdot E_{k}\leq 0, $ for all $k=1,...,n$, then $E$ has full support,
 i.e. $n_k>0$ for all $k=1,...,n$.
 \item If $\displaystyle E=\sum_{k=1}^{n} n_kE_{k}$ with 
 $n_k\in \bn^* $ for $k=1,...,n$, is an exceptional divisor such that 
$E\cdot E_{k} \leq -2K_{\tilde X}\cdot E_{k}, $ then for any $\alpha \in \bn^*$ we have 
$(\alpha E)\cdot E_{k} \leq -2K_{\tilde X}\cdot E_{k}. $
\end{itemize}

\end{remark}
\begin{definition}
Let $(X,O)$ be a germ of normal surface singularity, $\pi:\tilde X \longrightarrow X$ be the minimal
resolution of singularities, $E_{1},\ldots,E_{n}$ be the components 
of the exceptional divisor and $A=(a_{ i,j} )$ with $a_{ i,j}=E_i\cdot E_j$, be the $n\times n$
 symmetrical intersection  matrix  of 
the exceptional set of $\tilde X$. The dual graph $\Gamma$ of the intersection matrix $A$ is defined as follows:
\begin{itemize}
\item The vertices of the graph $\Gamma$ are  $E_{1}$, $\ldots,$  $E_{n}$,
\item For  $i\not= j$ there is an edge between $E_{i},$ and $E_{j}$ if and only if $a_{ i,j} \not= 0.$
\end{itemize}

\begin{remark} The graph $\Gamma$ is connected and conversely by a theorem due to Grauert, 
given  a 
 $n\times n$ symmetrical negative definite  matrix  $A=(a_{ i,j} )$ with a connected graph  
there exist a singularity  with $A$ as intersection  matrix.
 \end{remark}

\end{definition}
Now we introduce the definition of Nash numerical conditions, this is 
the central point of this work, in the other sections we will prove that
 Nash numerical conditions depend only on the  intersection matrix of
 the exceptional set.  A Nash matrix will be a matrix 
satisfying the Nash numerical conditions.  In  section 2, 3 
we characterize some Nash matrix, in section 4 we consider like star shaped
 graphs and in section 5  we present some examples.
 
\begin{definition}
Let $(X,O)$ be a germ of normal surface singularity, $\pi:\tilde X \longrightarrow X$ be the minimal
resolution of singularities and $E_{1},\ldots,E_{n}$ be the components 
of the exceptional divisor. Let $K_{\tilde X}$ the canonical divisor 
on ${\tilde X}$. We say that $(X,O)$ satisfies numerical Nash condition for $(i,j)$ if the following condition is fulfilled
$$\leqno (NN_{(i,j)}) \hskip 1cm\exists  E=\sum_{k=1}^{n} n_kE_{k}
 , n_k\in \bn^*  {\rm \ with\ } n_{i}<n_{j}  {\rm \ and \ } 
-E\cdot E_{k} \geq 2K_{\tilde X}\cdot E_{k}, \forall    k=1,\ldots,n$$
We also say that $(X,O)$ satisfies numerical Nash  condition, (NN), if $(NN_{(i,j)})$ 
is true for all couples $(i,j)$, with $i\not= j$. 
\end{definition}

As an immediate consequence of Proposition 1  and Theorem 1  we have:

\begin{corollary} With the above notations,   if $(X,O)$ satisfy numerical Nash condition 
for $(i,j)$
 then $ \bar{ \cal 
N}_{i}\not\subset \bar{ \cal N}_{j} $. In particular if $(NN$) is true then the Nash problem 
on arcs has a positive answer.
\end{corollary}
\begin{proposition}With the notations as above. Let $\Gamma$ be the dual 
graph of the intersection matrix of the exceptional set.
If $(NN$) is true for  $\Gamma$,  then
\begin{itemize}
\item  $(NN$) is true for any subgraph of  $\Gamma$
\item  $(NN$) is true by decreasing the self intersection numbers.
\end{itemize}
 \end{proposition}
\demo \begin{itemize}
\item  Let consider a subgraph $G$ of  $\Gamma$ and let $I$ be its support. Since $(NN)$ is true for  $\Gamma$, for any $i,j\in I, i\not= j$,  there exist 
$E=\sum_{k=1}^{n} n_kE_{k}
 , n_k\in \bn^*  {\rm \ with\ } n_{i}<n_{j}$ such that $ 
E\cdot E_{k} \leq -2K_{\tilde X}\cdot E_{k}, \forall    k=1,\ldots,n$

It then follows that for any $k\in I$,  $$ 
(\sum_{l\in I}n_lE_{l})\cdot E_{k} \leq -2K_{\tilde X}\cdot E_{k}-\sum_{l\notin I}
 n_lE_{l}\cdot E_{k} \leq -2K_{X'}\cdot E_{k},$$
where $K_{X'}$ is the canonical divisor 
of the minimal resolution singularity $X'$, having $G$ as  dual graph of the exceptional set. Remark that 
$K_{X'}\cdot E_{k}=K_{\tilde X}\cdot E_{k} $.
\item  In order to prove the second assertion it will be enough to consider one index $k\in \{1,...,n\}$ 
and the  intersection matrix $A'=(a'_{i,j})$
 defined by $a'_{i,j}=a_{i,j}$ if $(i,j)\not=(k,k)$ and $a'_{k,k}=a_{k,k}-1$. Let remark that the matrix $A'$ corresponds to 
 a minimal resolution of some isolated singularity, $\pi' : \tilde 
X'\longrightarrow X'$, call $E'_{1},\ldots, E'_{n}$ the irreducible components of the exceptional set in $\tilde 
X'$ (In fact as a curve $E'_{i}=E_{i}$, but we need to distinguish them in $\tilde X$ and $\tilde X'$. Let 
$E=\sum_{k=1}^{n} n_kE_{k}
 , n_k\in \bn^*  {\rm \ with\ } n_{i}<n_{j}$ such that $ 
E\cdot E_{k} \leq -2K_{\tilde X}\cdot E_{k}, \forall    k=1,\ldots,n$ and set $E'=\sum_{k=1}^{n} n_kE'_{k}$. By the Remark 1
we can assume that $n_{k}\geq 2$ for any $k=1,\ldots,n$.
It follows that

$$\begin{array}{ccl}
 K_{\tilde X'}\cdot E'_{i}&=&K_{\tilde X}\cdot E_{i}\ \ {\rm for \ }i\not= k  \\
K_{\tilde X'}\cdot E'_{k}&=&K_{\tilde X}\cdot E_{k}+1 \\
E'\cdot E'_{i}&=&E\cdot E_{i}\leq -2K_{\tilde X}\cdot E_{i}=-2K_{\tilde X}\cdot E_{i} \ \ {\rm for \ }i\not= k \\
E'\cdot E'_{k}&=&E\cdot E_{k}-n_k\leq -2K_{\tilde X}\cdot E_{k}-n_k= -2K_{\tilde X'}
\cdot E'_{k}-n_k+2\leq -2K_{\tilde X'}\cdot E'_{k}\\
\end{array}$$
This complete the proof of the second assertion.
\end{itemize}

\section{Nash matrices, Gauss sequences}
Let $\pi : \tilde 
X\longrightarrow X$ be the minimal resolution of singularities and let 
$A=(a_{ i,j} )$ be the $n\times n$ symmetrical intersection  matrix  of 
the exceptional set of $X$, consider an  exceptional effective divisor 
$E= x_{1}E_{1}+\ldots+x_{n}E_{n}$,then 
 $$E\cdot E_{k}= x_{1}E_{1}\cdot E_{k}+\ldots+x_{n}E_{n}\cdot E_{k}=
 x_{1}a_{k,1}+\ldots+x_{n}a_{k,n}.$$ Set  $^tX=(x_{1}¥,\ldots,x_{n})$ 
and $ ^tC=(-2K_{\tilde X}\cdot 
 E_{1},\ldots,- 2K_{\tilde X}\cdot E_{n}) and c_i=- 2K_{\tilde X}\cdot E_{i}$,  then
 \begin{enumerate}
 \item Corollary 1 can be translated into linear algebra:

If  the inequality: $A X\leq C$ has a solution $(x_{1}¥,\ldots,x_{n})\in \bn^n$ 
such that   $x_i< x_j$, then  $ \bar{ \cal 
N}_{i}\not\subset \bar{ \cal N}_{j} $
\item  The condition  $(NN_{(i,j)})$ is equivalent to the condition:

 the inequality : $A X\leq C$ has solutions $(x_{1}¥,\ldots,x_{n})\in \bn^n$ 
such that  $x_i< x_j$.
 \end{enumerate}

 Remark that since $\tilde X$ 
 is the minimal resolution we have $K_{\tilde X}\cdot E_{i}\geq 0$ for 
 any $i$.
In what  follows we allow the  intersection matrix $A$ to have rational terms, remark 
that after multiplication by a convenient integer it will correspond to a singularity.
\begin{lemma}
        Let $(X,O)$ be a germ of a normal surface singularity ,  $\pi : \tilde X\longrightarrow X$ 
be the minimal 
resolution of singularities. Assume that $\pi$ has only two exceptional components 
 $E_{1},E_{2}$. Let $A=\pmatrix{ -a &c \cr c &-b }$ the intersection matrix of $E_{1},E_{2}$. Then 
\begin{enumerate}
\item   $c<a$ if and only if $(NN_{(1,2)})$ is true
\item   $c<b$ if and only if $(NN_{(2,1)})$ is true
\item   $c<\min \{a,b\}$ if and only if $(NN)$ is true.
\end{enumerate}
In particular since   the quadratic form associated to the  matrix $A$ is negative definite, we have $c^2<ab$, 
which implies that either $ \bar{ \cal 
N}_{1}\not\subset \bar{ \cal N}_{2} $or $ \bar{ \cal 
N}_{2}\not\subset \bar{ \cal N}_{1} $.
\end{lemma}
\demo We are looking for solutions  $(x,y)\in \bn^*$ of the system:
\begin{center}
$-a x+cy\leq c_{1}\leq 0$

 (*) \hskip 4cm
 \phantom{$-a x+cy\leq c_{1}\leq 0$}
 
$cx-by\leq c_{2}\leq 0$
\end{center}
let $D_1$  the line of equation $-ax+cy=c_{1}$ and $D_2$ the line with equation $cx-by=c_{2}$, 
since $A$ is negative definite we have   $c^2<ab$, which implies  $c/b<a/c$, so the relative positions of the  lines $D_1,D_2$, 
and the set of solutions of the system (*) are  represented in figures below. Since 
these are the unique possible cases we are done.

\includegraphics[height=1.8in]{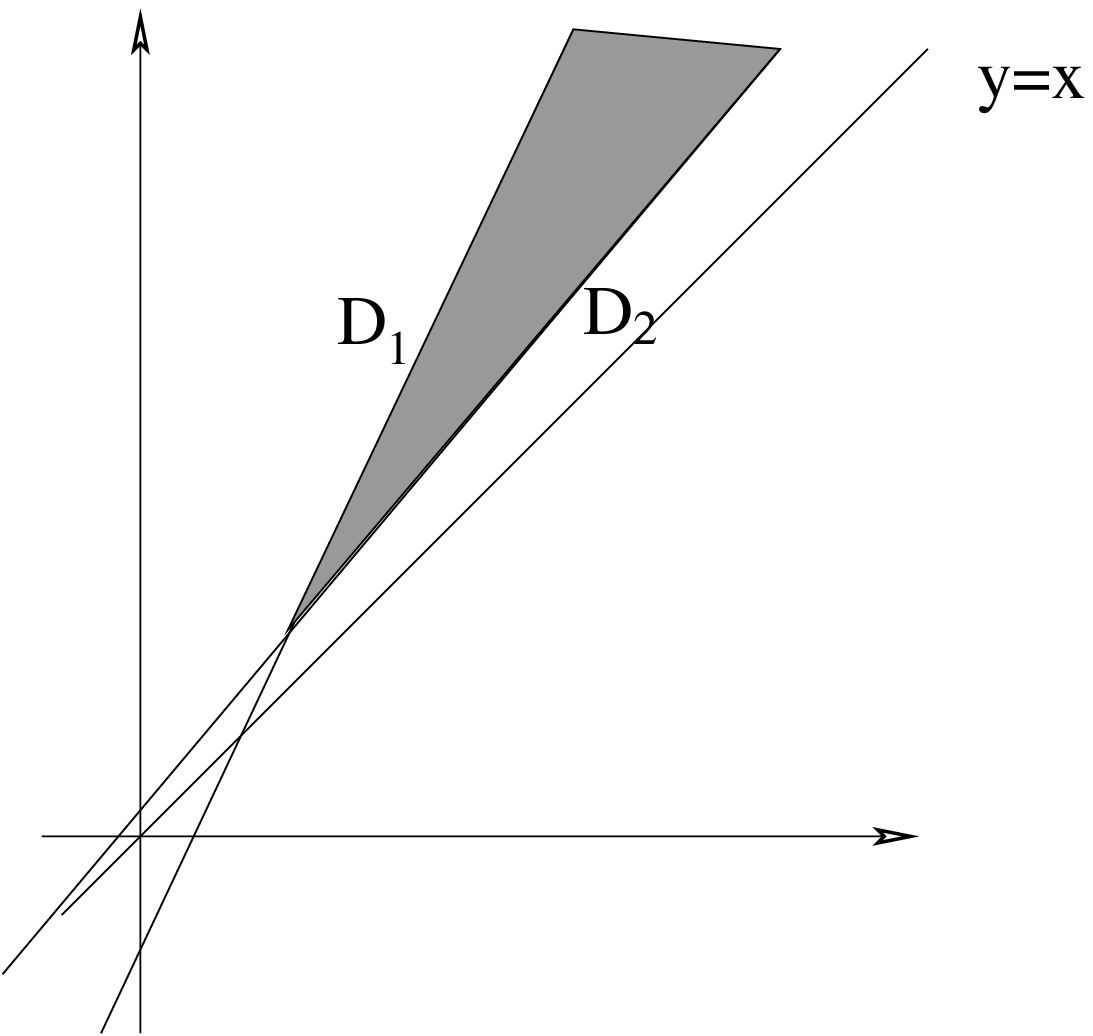}
\includegraphics[height=1.8in]{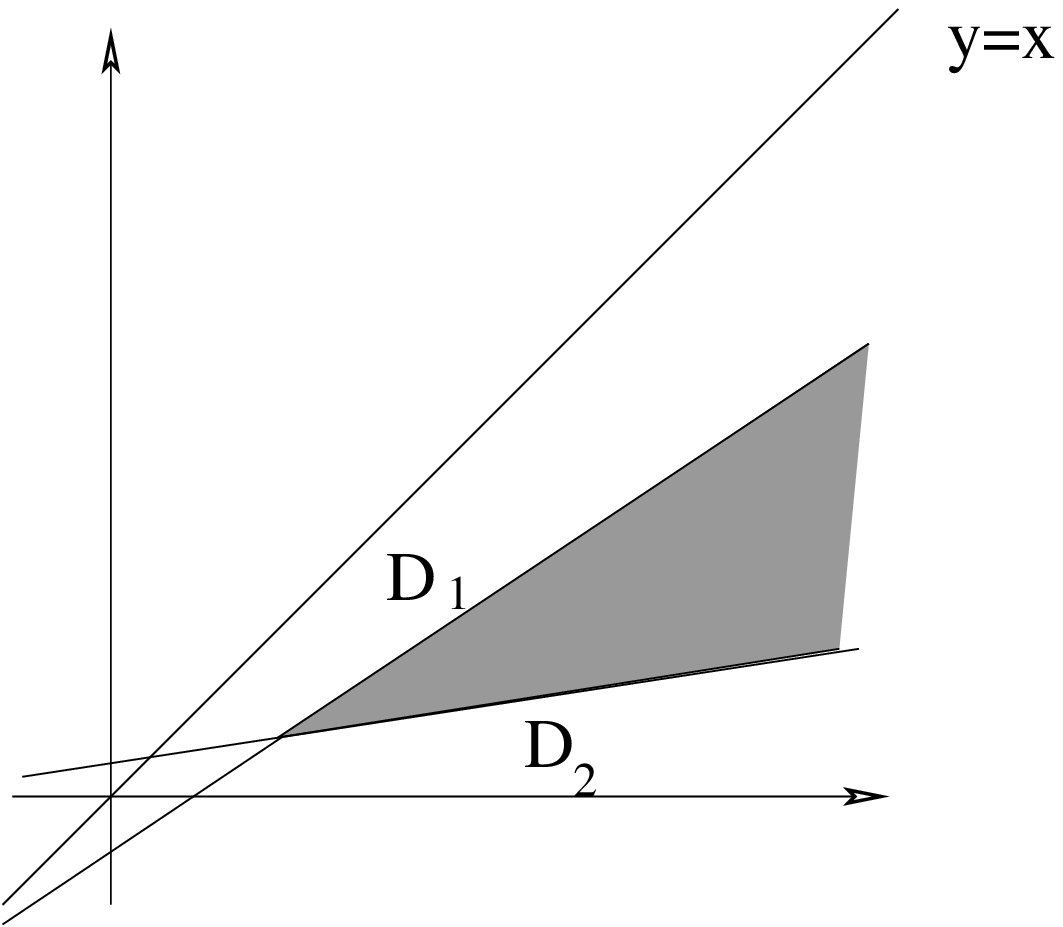}      
\includegraphics[height=1.8in]{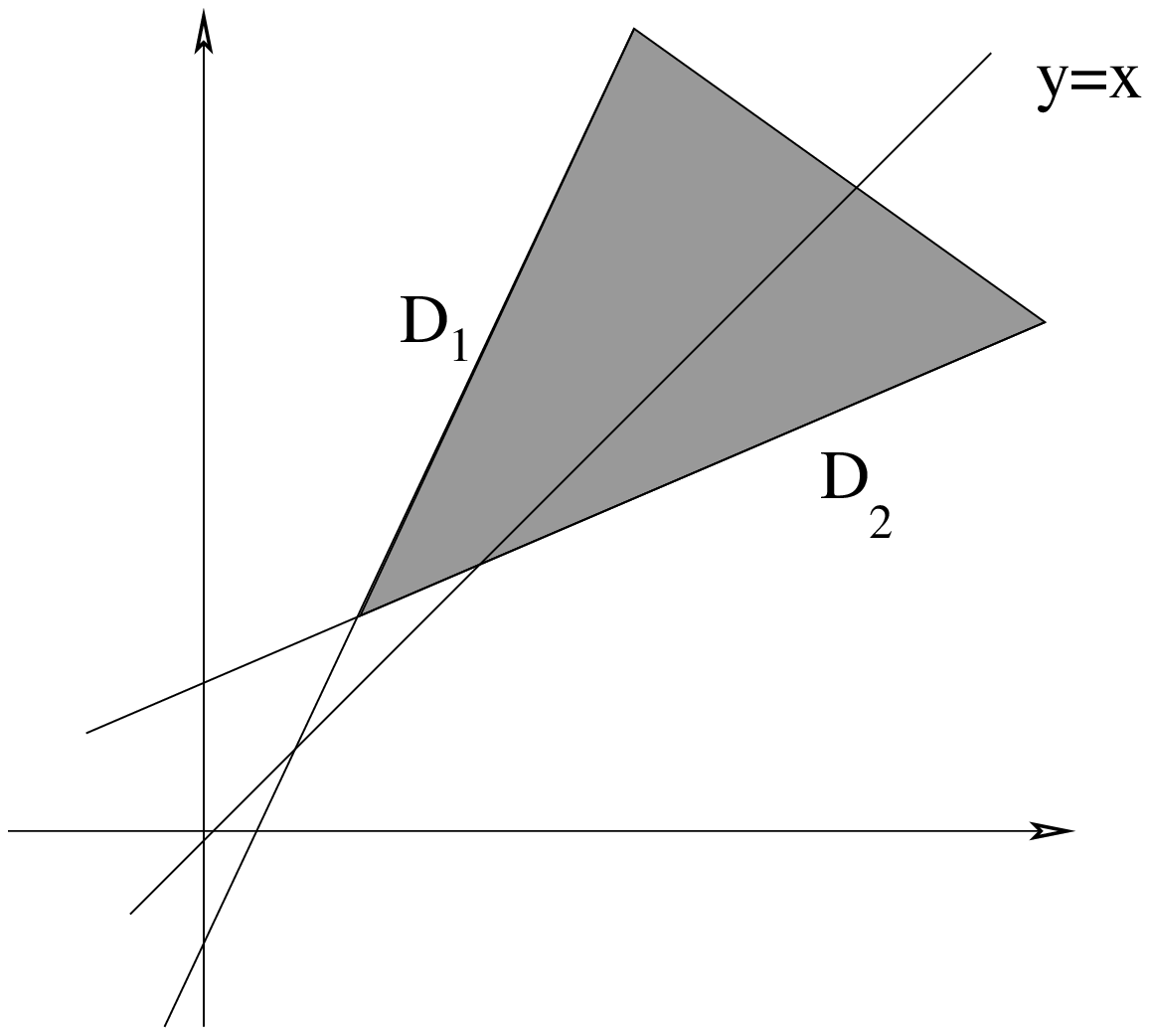}  
 $\phantom{..} \hskip 0.5cm (i))\  c<a$ and  $b\leq c$\hskip 2.3cm $(ii)\ c<b$ but $a\leq c$\hskip 2.5cm
  $(iii)\ c<a$ and  $c<b$.   
\begin{corollary}
Let  $\pi 
: \tilde X\longrightarrow X$ be the minimal resolution of 
singularities and $E_{1},\ldots,E_{n}$ be the components of the 
exceptional divisor,then for any $i\not = j$ either $ \bar{ \cal 
N}_{i}\not\subset \bar{ \cal N}_{j} $or $ \bar{ \cal 
N}_{j}\not\subset \bar{ \cal N}_{i} $. In any case if $i\not = j$ then  $ \bar{ \cal 
N}_{i}\not= \bar{ \cal N}_{j} $. In particular after considering numerical Nash conditions, 
in order to check if Nash is true, we will be 
reduced to check at most the half of non inclusion conditions.
\end{corollary}

We prove the Corollary by induction on $n$. For $ n=2$ it was proved in lemma 1.

\noindent Assume $n\geq 3$, by changing the order in the set $E_{1},\ldots,E_{n}$, 
we can suppose that $i=1$ and $j=2$, now    pick  $k$ a positive integer such that $ka_{n,n}<c_n$ and put   
$\displaystyle -a_{n,n}x_{n}=\sum_{1\leq i\leq n-1}a_{n,i} x_{i}-
ka_{n,n} $  in our system $AX\leq C$,  then we have 
the inequality: $A' X \leq C'$ where $
a'_{i,j}=a_{n,i}a_{n,j} -a_{i,j}a_{n,n}$ for all $i, j$ and $c'_{j}= (-c_{j}+ka_{n,j})a_{n,n}$.
By induction hypothesis there exist a vector 
 $S=(s_{1},\ldots, s_{n-1})\in \bn^n$  solution of the in-equation $$A' X\leq 
 C'$$ with $s_{1}\not= s_{2}$.  Let 
 $s_{n}=a_{n,1}s_{1}+\ldots+a_{n,n-1}s_{n-1}+k$, 
 then a simple computation shows that  the vector $T=(-a_{n,n}s_{1},\ldots, -a_{n,n}s_{n-1},s_{n})$  is a  solution of $A X \leq C$ for $k$ large enough.
 
 \noindent Remark that by construction the vector $T$ has strict 
 positive components.

Now we consider the sequences appearing in the proof of the last Corollary.
\begin{definition} Let:
 $a^{(n)}_{i,j}=a_{i,j}$ and for  any
 $2\leq l \leq n-1  $ set 
$\displaystyle a^{(l)}_{i,j}= a^{(l+1)}_{i,j} - 
{{a^{(l+1)}_{l+1,i}a^{(l+1)}_{l+1,j}}\over{a^{(l+1)}_{l+1,l+1}} } 
,\quad  1\leq i,j\leq  l $. Also for  any
 $2\leq l \leq n $ let $ C(A)_i^{(l)}=\sum_{j=1}^l a^{(l)}_{ i,j}$. 
 We will also use the notation $ C(A)_i = C(A)_i^{(n)}.$
\end{definition} 

\begin{lemma}The matrices  $A^{(l)}=(a^{(l)}_{ i,j} )$ appear naturally when
 we use the Gauss method to decompose the quadratic form associated to $A$ 
into a  sum of squares. In particular the matrix $A^{(l)}$ are negative definite.
For this reason we will call the terms $a^{(l)}_{ i,j}$ 
the Gauss sequence associated to $A$. 
\end{lemma}
\demo The quadratic form associated to the matrix $A$ is:
$$Q= \sum_{i=1}^n a_{ i,i}x_i^2 +2 \sum_{1\leq i<j\leq n} a_{ i,j}x_ix_j$$ we follow Gauss method to 
squaring a quadratic form:
 $$Q= \sum_{i=1}^{n-1} a_{ i,i}x_i^2 +2 \sum_{1\leq i<j\leq n-1}a_{ i,j}x_ix_j +
a_{ n,n}x_n^2 +2 \sum_{i=1}^{n-1} a_{ i,n}x_ix_n $$
but $$ a_{ n,n}x_n^2 +2 \sum_{i=1}^{n-1} a_{ i,n}x_ix_n = 
a_{ n,n}(x_n+\sum_{i=1}^{n-1} {{a_{ i,n}}\over{a_{ n,n}}}x_i)^2-
\sum_{i=1}^{n-1}{{a_{ i,n}^2}\over{a_{ n,n}}}x_i^2-
2\sum_{1\leq i<j\leq n-1}{{a_{ i,n}a_{ j,n}}\over{a_{ n,n}}}x_ix_j
 $$
 Hence $$ Q= a_{ n,n}(x_n+\sum_{i=1}^{n-1} {{a_{ i,n}}\over{a_{ n,n}}}x_i)^2+ 
\sum_{i=1}^{n-1} (a_{ i,i}- {{a_{ i,n}^2}\over{a_{ n,n}}})x_i^2 +
2 \sum_{1\leq i<j\leq n-1} (a_{ i,j}-{{a_{ i,n}a_{ j,n}}\over{a_{ n,n}}})x_ix_j $$
and $A$ is negative definite if and only if $a_{ n,n}<0$ and  $A^{(n-1)}$ 
is negative definite.
\begin{remark}
  \begin{enumerate}
  
\item By multiplying by a convenient natural number the matrix $A$ has integer coefficients and correspond to some singularities. Our definition 
does not depend on the topological type of the components of the exceptional divisor.
\item For $l\geq 3$ the operation  $A^{(l)}\mapsto  A^{(l-1)}$ consist to contract the exceptional component $E_l$ in the graph $\Gamma_l$ corresponding to 
 $A^{(l)}$, it is an algebraic operation and this contraction has no geometry meaning. In what follows we will use this notation. 
  \end{enumerate}
\end{remark}

We have immediately from lemma 1 and Corollary 2 that
\begin{proposition}
        
Let $\pi : \tilde 
X\longrightarrow X$ be the minimal resolution of singularities and let 
$A=(a_{ i,j} )$ be the $n\times n$ symmetrical intersection  matrix  of 
the exceptional set of $X$. Then 
\begin{enumerate}
\item   $a^{(2)}_{1,2}<-a^{(2)}_{1,1} $ if and only if $(NN_{(1,2)})$ is true
\item   $a^{(2)}_{1,2}<-a^{(2)}_{2,2} $ if and only if $(NN_{(2,1)})$ is true
\item   $a^{(2)}_{1,2}<\min \{-a^{(2)}_{1,1},-a^{(2)}_{2,2} \}$ if and only if both 
$(NN_{(1,2)}), (NN_{(2,1)})$  are true.
\end{enumerate}

\end{proposition}
\begin{theorem} Let $\pi : \tilde 
X\longrightarrow X$ be the minimal resolution of singularities, let 
$A=(a_{ i,j} )$ be the $n\times n$ symmetrical intersection  matrix  of 
the exceptional set of $X$ and let $ C(A)_i^{(l)}=\sum_{j=1}^l a^{(l)}_{ i,j}$. For $l\geq 1$, we 
consider the property:
$$\leqno{(*_{l+1})} C(A)_i^{(l+1)}<0, {\rm \ for\ \ }   i=1,...,l+1.$$ 
If $(*_{l+1})$ is true for  some $l\geq 2$ then $(*_{l})$ is true.

Let $\sigma\in S_n$ any permutation of $E_1,...,E_n$, we denote by  $A^\sigma$ 
the corresponding 
intersection  matrix 
 obtained from $A$ by permuting lines and columns. Then 
$(N N)$ is true if and only if there exist a natural integer $l\geq 1$ such that 
$$\leqno{(*_{l+1})} C(A^\sigma)_i^{(l+1)}<0, {\rm \ for\ \ }   i=1,...,l+1, \forall \sigma \in S_n.$$ 
   In particular we recover the following result from  \cite{P-PP}: if $ C(A)_i^{(n)}<0, {\rm \ for\ \ }   i=1,...,n $ then the Nash map 
${ \cal N }$ is bijective.

Note that condition $(*_{l})$  has a meaning only if $l\geq 2$.

\end{theorem}

\demo Assume that $  C(A)_i^{(l+1)}<0, {\rm \ for\ \ }   i=1,...,l +1$, let $i\leq l$, 
by definition $$C(A)_i^{(l)}=\sum_{j=1}^l a^{(l)}_{ i,j}=\sum_{j=1}^{l}
(a^{(l+1)}_{i,j} - 
{{a^{(l+1)}_{l+1,i}a^{(l+1)}_{l+1,j}}\over{a^{(l+1)}_{l+1,l+1}} } )
$$
$$C(A)_i^{(l)}=\sum_{j=1}^{l}a^{(l+1)}_{i,j}-{{a^{(l+1)}_{l+1,i}}\over{a^{(l+1)}_{l+1,l+1}} } 
\sum_{j=1}^{l}a^{(l+1)}_{l+1,j}
,
$$
$$C(A)_i^{(l)}= C(A)_i^{(l+1)}- {{a^{(l+1)}_{l+1,i}}\over{a^{(l+1)}_{l+1,l+1}} } 
C(A)_{l+1}^{(l+1)}  <0
$$
The second assertion follows from Proposition 3.
Remark that it is not necessary to consider all permutation of $E_1,...,E_n$.
 
\begin{definition} Let 
$A=(a_{ i,j} )$ be the $n\times n$ symmetrical negative definite matrix with rational 
coefficients with $a_{ i,i}<0, a_{ i,j}\geq 0$ for all $i,j, i\not= j$. 
We say that $A$ is a Nash matrix if for any permutation $\sigma$ of the set 
$\{1,...,n\}$ $C(A^\sigma)_{1}^{(2)}<0,C(A^\sigma)_{2}^{(2)}<0 $ 
  
\end{definition}

\section{Trees, Cycles, Generalized Cycles}
We look now for some necessary or sufficient conditions   in order to have
 the condition $(N N)$ true. For the moment we need to recall some notation on graphs.
 
 \begin{definition} Let 
$A=(a_{ i,j} )$ be the $n\times n$ symmetrical negative definite matrix with rational 
coefficients with $a_{ i,i}<0, a_{ i,j}\geq 0$ for all $i,j, i\not= j$. 
Let $\Gamma$ the dual graph associated to $A$. 
We say that $E_i$ is a leaf of $\Gamma$ 
if $a_{ i,j}\not= 0$ for exactly one index $j\not= i$ i.e. $E_i$ is connected 
to only one other vertex of $\Gamma$. A cycle of $\Gamma$ is a  subgraph  $ {\cal C}$ 
 where every vertex is
 connected to exactly  two others  vertex of $ {\cal C}$. A tree is a subgraph with no cycles. 
Finally a complete subgraph is a subset of $\Gamma$, where every two points are connected.

\end{definition}

\begin{lemma} Assume that for any point $E_j$ of  $\Gamma$,we have $ C(A)_j \leq 0$.
\begin{enumerate}
\item  For any $l\leq n$ and $j\leq l$ we have  $ C(A)^{(l)}_{j} \leq 0$
\item  If $ C(A)^{(l+1)}_{i} < 0$   then  $ C(A)^{(l)}_{i} < 0$
\item Let consider a path $E_{i_1},E_{i_2},...,E_{i_k}$ in $\Gamma$, and   $ C(A)_{i_k} < 0$. 
After contracting $E_{i_k},E_{i_{k-1}},...,E_{i_2}$we will have $ C(A)^{(2)}_{i_1} < 0$.

 \end{enumerate}
\end{lemma}
\demo The first two assertions follow  immediately from the following formula, 
which is true for  any $l\geq 2$, and $1\leq i\leq l$:
$$C(A)_i^{(l)}= C(A)_i^{(l+1)}- {{a^{(l+1)}_{l+1,i}}\over{a^{(l+1)}_{l+1,l+1}} } 
C(A)_{l+1}^{(l+1)}  <0
$$
We prove the third assertion  by induction on $k$ the length of the path, if $k=2$, 
by the above formula we get the answer. Now take any $k\geq 3$,
then using again the above formula we have that $ C(A)^{(n-1)}_{i_{k-1}} < 0$, 
 by the induction hypothesis we get $ C(A)^{(n-k+1)}_{i_{1}} < 0$, 
so by the assertion 1 we  are done.
\begin{theorem}Let $\pi : \tilde 
X\longrightarrow X$ be the minimal resolution of singularities and let 
$A=(a_{ i,j} )$ be the $n\times n$ symmetrical intersection  matrix  of 
the exceptional set of $X$. If $(N N)$ is true then  $ C(A)_i < 0$ for any leaf 
$E_i$ of  $\Gamma$. 

\end{theorem}
\demo Suppose that $(N N)$ is true.  Let  $E_i$ be a leaf of  $\Gamma$, 
we can assume that $i=1$ and 
   $E_2$ is the unique   vertex connected to $E_1$, by contracting all other vertex of $\Gamma$,
 we will have $ a^{(2)}_{1,1}= a_{1,1}, a^{(2)}_{1,2}= a_{1,2} $. By Proposition 3 (or Theorem 2)  
we must have 
$ C(A)_1= C(A)^{(2)}_1< 0$. This concludes the proof.

\begin{theorem} Let $\pi : \tilde 
X\longrightarrow X$ be the minimal resolution of singularities and let 
$A=(a_{ i,j} )$ be the $n\times n$ symmetrical intersection  matrix  of 
the exceptional set of $X$. Assume that $\Gamma$ is a tree and 
$ C(A)_i \leq 0$ for any vertex $E_i$ of $\Gamma$.
 Then  
$ C(A)_i < 0$ for any leaf  $E_i$ of $\Gamma$
if and only if 
$(N N)$ is true. In particular if the above conditions are satisfied the Nash map 
${ \cal N }$ is bijective.

\end{theorem}
\demo The necessary condition was proved before.  
The proof of the other implication is by induction on $n$. 
If $n=2$ the hypothesis implies that $(N N)$ is true by Lemma 1.
So assume the case $n-1$ is solved and we prove the case $n$. Take any $i\not= j$. 
We have two cases

1) Both $E_i, E_j$ are leaves of $\Gamma$, then $ C(A)^{(n)}_i<0,C(A)^{(n)}_j<0 $, 
by contracting all vertex in $\Gamma$ except $E_i,E_j$ and applying Lemma 2, we get that 
$ C(A)^{(2)}_i<0,C(A)^{(2)}_j<0 $, and we are done.

2) At most  one of $E_i,E_j$ is   a leaf, then there exist a leaf $E_k$, different 
 from $E_i$, $E_j$, so after changing the 
order of  the exceptional 
components we can assume that $i=1,j=2,k=n$, let $E_l$ be unique component connected to $E_n$.
By contracting $E_n$, we get the matrix $A^{(n-1)}=(a^{(n-1)}_{i,j})$, 
with $a^{(n-1)}_{i,j}= a_{i,j}$ for any $(i,j)\not= (l,l)$ and 
$\displaystyle a^{(n-1)}_{l,l}=a_{l,l}-{{a_{l,n}^2}\over{a_{n,n}}}<0$. It follows that
 $ C(A)^{(n-1)}_i=C(A)_i \leq 0$ for $i\not= l$ and $\displaystyle  C(A)^{(n-1)}_l=  
C(A)_l-({{a_{l,n}^2}\over{a_{n,n}}}+a_{l,n})<C(A)_l\leq 0$. Also the graph 
corresponding to the matrix  $A^{(n-1)}$ is a tree, so by induction hypothesis 
$(NN_{(1,2)}) $  and $(NN_{(2,1)}) $  are true, and we are done.
\begin{remark} Inside the class of rational singularities, rational minimal singularities are exactly those for which the graph 
satisfies the hypothesis of the above theorem. Note that Nash problem's on arcs for (rational) 
minimal  
singularities has a positive solution by the work 
of  Ana Reguera \cite{R}, also C. Plenat \cite{P1} and Fernandez-Sanchez \cite{F} gave different proofs. Our Theorem
 applies without any restriction 
on the topological type of the exceptional components and so extends to non rational singularities the mentioned results. 
\end{remark}
\begin{theorem} Let $\pi : \tilde 
X\longrightarrow X$ be the minimal resolution of singularities and let 
$A=(a_{ i,j} )$ be the $n\times n$ symmetrical intersection  matrix  of 
the exceptional set of $X$. Assume that  the graph $\Gamma$ of the exceptional set is a cycle,
 with $n\geq 3$,  and  $ C(A)_i \leq 0$ for all $i$.
  Then 
$(N N)$ is true if and only if $ C(A)_i < 0$ for at least two exceptional components.
In particular if these conditions are fulfilled the Nash map 
${ \cal N }$ is bijective.
\includegraphics[height=1in]{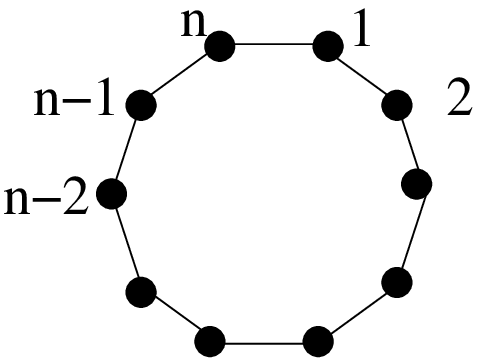}
\end{theorem}
\demo Assume first that  $ C(A)_i < 0$ for at least two exceptional components
We prove that $(N N)$ is true by induction on $n$.

If $n=3$, we contract the exceptional fiber $E_3$ and we get the matrix:
$$A^{(2)}= \pmatrix{a_{1,1}-{{a_{1,3}^2}\over{a_{3,3}}} &
 a_{1,2}-{{a_{1,3}a_{2,3}}\over{a_{3,3}}}\cr
  a_{1,2}-{{a_{1,3}a_{2,3}}\over{a_{3,3}}}&
a_{2,2}-{{a_{2,3}^2}\over{a_{3,3}}}\cr }$$

 It follows that $\displaystyle  C(A)^{(2)}_{1}=  
C(A)_{1}-({{a_{1,3}}\over{a_{3,3}}}(C(A)_{3}))<  0$ and
$\displaystyle  C(A)^{(2)}_{2}=  
C(A)_{2}-({{a_{2,3}}\over{a_{3,3}}}(C(A)_{3}))<  0$
 since by hypothesis two over the three numbers
 $C(A)_{1},C(A)_{2},C(A)_{3}$ are strictly negative. So the case $n=3$ is over.

Consider now the case $n\geq 4$.  By contracting   $E_n$,
 we get the  matrix  
 $A^{(n-1)}=(a^{(n-1)}_{i,j})$, 
with $a^{(n-1)}_{i,j}= a_{i,j}$ if $i,j\in \{ 2,3,..., n-2\}$ and 
$\displaystyle a^{(n-1)}_{1,1}=a_{1,1}-{{a_{1,n}^2}\over{a_{n,n}}}$, 
$\displaystyle a^{(n-1)}_{1,n-1}=-{{a_{1,n}a_{n-1,n}}\over{a_{n,n}}}$ and
$\displaystyle a^{(n-1)}_{n-1,n-1}=a_{n-1,n-1}-{{a_{n-1,n}^2}\over{a_{n,n}}}$.

 It follows that
 $ C(A)^{(n-1)}_i=C(A)_i \leq 0$ for $i\in \{ 2,3,..., n-2\}$, 
 $ \displaystyle C(A)^{(n-1)}_{n-1}=  
C(A)_{n-1}-({{a_{n-1,n}}\over{a_{n,n}}}(C(A)_{n}))$ and
  $ \displaystyle C(A)^{(n-1)}_{1}=  
C(A)_{1}-({{a_{1,n}}\over{a_{n,n}}}(C(A)_{n})) $    
                         We have to consider three cases:

i)  $ C(A)_{1}=C(A)_{n-1}=C(A)_{n}=0$ then  there are two indexes $i,j\in \{ 2,3,..., n-2\}$ such that
 $  C(A)^{(n-1)}_i<0, C(A)^{(n-1)}_j<0$.

ii)  $C(A)_{n}<0$, 
then $C(A)^{(n-1)}_{n-1}, C(A)^{(n-1)}_{1}$ are strictly negative. 
         
iii)  at least one of  $ C(A)_{1}=0$ and $C(A)_{n-1} $   and $C(A)_{n}=0$, 
then either  $C(A)^{(n-1)}_{n-1}<0$ or $ C(A)^{(n-1)}_{1}<0$.
         
So the induction hypothesis is verified by $A^{(n-1)}$ and   we are done.

Conversely, if $(N N)$ is true and   $ C(A)_i < 0$ for at most one index $i$, 
take any index $j\not= i$, by contracting all other components $E_k, k\not= i,j$ we will have 
$C(A)^{(2)}_i=C(A)_i=0, C(A)^{(2)}_j=C(A)_j<0 $, this is a contradiction by Proposition 3.
We can give a more general result that the preceding one, for this we need some definitions. 
\begin{definition}
We say that a subgraph $G$ of $\Gamma$ is a generalized cycle if any two vertex
 of $G$ are connected by a cycle. Remark that a cycle or a complete graph are 
generalized cycles.

A generalized cycle is a leaf of  $\Gamma$ if at most one vertex of $G$ 
is connected to one vertex of $\Gamma \setminus G$.
\end{definition}
The proof of the next Corollary is exactly the same as for a cycle, 
and we left it to the reader:
\begin{corollary} Let $\pi : \tilde 
X\longrightarrow X$ be the minimal resolution of singularities of $X$ and let 
$A=(a_{ i,j} )$ be the $n\times n$ symmetrical intersection  matrix  of 
the exceptional set of $X$. Suppose that $\Gamma$ is a generalized cycle.

We assume that $n\geq 3$ and 
$C(A)_i\leq 0, $ for any vertex $E_i$. Then 
$(N N)$ is true if and only if $C(A)_i<0, $ for at least 
two vertex.

\end{corollary}
\begin{example}The following matrix and  graph correspond to a generalized cycle,  for which Nash's problem has an affirmative answer.
\vskip 0.5cm
 $ A=\pmatrix{-5 & 1 & 1 & 1 & 1 & 1 & 0 &
 0\cr 1 &
 -6 & 1 & 1 & 0 & 1 & 1 & 0\cr 
1 & 1 & -5 & 1 & 0 & 0 & 1 & 1\cr 1 & 1 & 1 & -5 & 1 & 0 & 0 & 1\cr
 1 & 0 & 0 & 1 & -2 & 0 & 0 & 0\cr 1 & 1 & 0 & 0 & 0 & -2 & 0 & 0\cr 
 0 & 1 & 1 & 0 & 0 & 0 & -2 & 0\cr
0 & 0 & 1 & 1 & 0 & 0 & 0 & -3}$
\setlength{\unitlength}{2144sp}%
\begingroup\makeatletter\ifx\SetFigFont\undefined%
\gdef\SetFigFont#1#2#3#4#5{%
  \reset@font\fontsize{#1}{#2pt}%
  \fontfamily{#3}\fontseries{#4}\fontshape{#5}%
  \selectfont}%
\fi\endgroup%
\begin{picture}(3495,0)(0,-1500)
{\thinlines
\put(4415,-1344){\circle*{134}}
}%
{\put(3706,-2034){\circle*{134}}
}%
{ \put(2984, 61){\circle*{134}}
}%
{ \put(3706,-616){\circle*{134}}
}%
{ \put(2266,-639){\circle*{134}}
}%
{ \put(1546,-1359){\circle*{134}}
}%
{ \put(2281,-2072){\circle*{134}}
}%
{ \put(2986,-2776){\circle*{134}}
}%
{ \put(4426,-1336){\line(-1, 1){1440}}
\put(2986,106){\line(-1,-1){1440}}
\put(1546,-1334){\line( 1,-1){1440}}
\put(2986,-2774){\line( 1, 1){1440}}
}%
{ \put(2311,-661){\line( 1,-1){1395}}
}%
{ \put(2266,-2056){\line( 1, 0){1440}}
\put(3706,-2056){\line( 0, 1){1440}}
\put(3706,-616){\line(-1, 0){1440}}
\put(2266,-616){\line( 0,-1){1440}}
}%
{ \put(2266,-2056){\line( 1, 1){1395}}
}%
\put(1540,-571){\makebox(0,0)[lb]{\smash{{\SetFigFont{10}{16.8}{\rmdefault}{\mddefault}{\updefault}{ -5}%
}}}}
\put(3830,-526){\makebox(0,0)[lb]{\smash{{\SetFigFont{10}{16.8}{\rmdefault}{\mddefault}{\updefault}{ -6}%
}}}}
\put(3920,-2191){\makebox(0,0)[lb]{\smash{{\SetFigFont{10}{16.8}{\rmdefault}{\mddefault}{\updefault}{ -5}%
}}}}
\put(1500,-2236){\makebox(0,0)[lb]{\smash{{\SetFigFont{10}{16.8}{\rmdefault}{\mddefault}{\updefault}{ -5}%
}}}}
\put(1000,-1426){\makebox(0,0)[lb]{\smash{{\SetFigFont{10}{16.8}{\rmdefault}{\mddefault}{\updefault}{ -2}%
}}}}
\put(4640,-1426){\makebox(0,0)[lb]{\smash{{\SetFigFont{10}{16.8}{\rmdefault}{\mddefault}{\updefault}{ -2}%
}}}}
\put(2800,239){\makebox(0,0)[lb]{\smash{{\SetFigFont{10}{16.8}{\rmdefault}{\mddefault}{\updefault}{ -2}%
}}}}
\put(2840,-3136){\makebox(0,0)[lb]{\smash{{\SetFigFont{10}{16.8}{\rmdefault}{\mddefault}{\updefault}{ -3}%
}}}}
\end{picture}%

\vskip 0.5 cm

\end{example}
\vskip 0.5cm
\begin{theorem} Let $\pi : \tilde 
X\longrightarrow X$ be the minimal resolution of singularities of $X$ and let 
$A=(a_{ i,j} )$ be the $n\times n$ symmetrical intersection  matrix  of 
the exceptional set of $X$.

We assume that $n\geq 3$, $\Gamma$ is not a generalized cycle and 
 \begin{enumerate}
\item $C(A)_i\leq 0, $ for any vertex $E_i$.
\item $C(A)_i< 0, $ for any leaf $E_i$ of $\Gamma$.
\end{enumerate}
Then 
$(N N)$ is true if and only if for any generalized cycle $G$ of  $\Gamma$ 
there is a vertex in $G$, not connected 
to one vertex of $\Gamma \setminus G$ such that $C(A)_i<0$.
\end{theorem}

\demo   
Assume that $(N N)$ is true, we have seen that  $C(A)_i< 0, $ for any leaf $E_i$ of $\Gamma$,
 now consider any 
generalized cycle $G$, we contract  all points outside this generalized cycle, 
so $(NN)$ is still true
 for this $G$, this implies that $C(A)_i<0$ for at least one vertex in $G$ not connected 
to one vertex of $\Gamma \setminus G$. We have finished to prove the necessary condition.
 
We prove now the other implication. Take two vertex $E_i , E_j$ of $\Gamma$, since $\Gamma$ is connected, 
there is a path ${\cal C}$  in $\Gamma$ connecting them. We must consider two cases, 
 
 \begin{enumerate}
\item ${\cal C}$ cannot  be extended   to a cycle, then by contracting all the vertex not
 in ${\cal C}$, we are reduced to the case of a tree, which was solved in Theorem 4.

\item ${\cal C}$ can be extended   to a cycle then $E_i , E_j$ are inside a generalized
 cycle, then by contracting all the vertex not
 in ${\cal C}$, we are reduced to the case of a generalized cycle, which was solved 
just before.
\end{enumerate}

\section{Like Star graphs}
We can improve the above  result in the some special situations:
\begin{theorem} Let $\pi : \tilde 
X\longrightarrow X$ be the minimal resolution of singularities and let 
$A=(a_{ i,j} )$ be the $n\times n$ symmetrical intersection  matrix  of 
the exceptional set of $X$. Assume that $X$ has a polygon singularity, i.e. 
 the graph of the exceptional set is a star
 with root $E_n$ and all other vertex are leaves.
\includegraphics[height=1in]{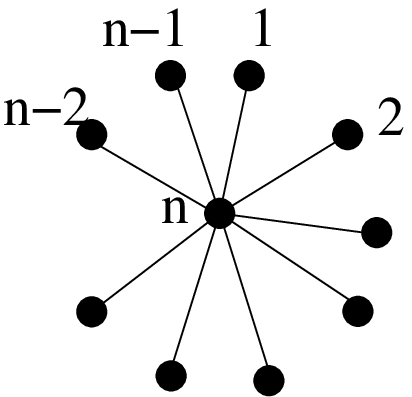} 

Then $(N N)$ is true if and only if we have the following conditions:
\begin{itemize}
\item $\forall i=1,...,n-1 \ \ \ a_{i,i}+a_{i,n}<0$
\item $\forall i,j=i,...,n-1\ \ \ a_{i,i}a_{j,j}\Delta_n + 
a_{j,n}(a_{i,i}a_{j,n}-a_{j,j}a_{i,n})<0$
\item $\displaystyle \forall i=i,...,n-1 \  \ \ {{a_{i,n}}\over{a_{i,i}}}(a_{i,i}+a_{i,n})+\Delta_n<0$

\end{itemize}
where $\displaystyle \Delta_n=a_{n,n}-\sum_{i=1}^{n-1}{{a_{i,n}^2}\over{a_{i,i}}}$. We note that $A$ is negative 
definite if and only if $\Delta_n<0$.
\end{theorem}
\demo
It is enough to compare any two leaves and any leaf with the root. So we consider the 
following order
 :$E_1,E_2, E_n,E_3,...,E_{n-1}$. After applying the construction above we are reduced to the matrix
 $$A^{(3)}= \pmatrix {a_{1,1}&0&a_{1,n}\cr
 0&a_{2,2}&a_{2,n}\cr
 a_{1,n}&a_{2,n}&a^{(3)}_{n,n}\cr}$$
 where $a^{(3)}_{n,n}=a_{n,n}-\sum_{i=3}^{n-1}{{a_{i,n}^2}\over{a_{i,i}}}$.
We are reduced to consider two cases:

first case: Comparison of $E_1,E_2$
Again by the construction above we are reduced to the matrix:
$$A^{(2)}= \pmatrix {a_{1,1}-{{a_{1,n}^2}\over{a^{(3)}_{n,n}}}&
-{{a_{1,n}a_{2,n}}\over{a^{(3)}_{n,n}}}\cr
-{{a_{1,n}a_{2,n}}\over{a^{(3)}_{n,n}}}&
a_{2,2}-{{a_{2,n}^2}\over{a^{(3)}_{n,n}}}\cr}
$$
So $(NN)_{(1,2)}$ and $(NN)_{(2,1)}$ are true if and only if :
$$a_{1,1}-{{a_{1,n}^2}\over{a^{(3)}_{n,n}}} -{{a_{1,n}a_{2,n}}\over{a^{(3)}_{n,n}}}<0$$
$$a_{2,2}-{{a_{2,n}^2}\over{a^{(3)}_{n,n}}}-{{a_{1,n}a_{2,n}}\over{a^{(3)}_{n,n}}}<0$$
by simple computations these are equivalent to:
$$a_{i,i}a_{j,j}\Delta_n + 
a_{j,n}(a_{i,i}a_{j,n}-a_{j,j}a_{i,n})<0$$
for $\{i,j\}=\{1,2\}$.

Let consider now the second case: Comparison of $E_1,E_n$

By the construction above we are reduced to the matrix:
$$ \pmatrix {a_{1,1}&a_{1,n}\cr
a_{1,n}&a^{(3)}_{n,n}-{{a_{2,n}^2}\over{a_{2,2}}}\cr}
$$
So $(NN)_{(1,2)}$ and $(NN)_{(2,1)}$ are true if and only if :
$a_{1,1}+a_{1,n}<0$
and
$a_{1,n}+a^{(3)}_{n,n}-{{a_{2,n}^2}\over{a_{2,2}}}<0$.

After simple computations these are equivalent to:
$a_{1,1}+a_{1,n}<0$
and
${{a_{1,n}}\over{a_{1,1}}}(a_{1,1}+a_{1,n})+\Delta_n<0$
Since the choice of the leaves were arbitrary, we are done.

The next  corollary follows immediately from the theorem.

\begin{corollary}Let $\pi : \tilde 
X\longrightarrow X$ be the minimal resolution of singularities and let 
$A=(a_{ i,j} )$ be the $n\times n$ symmetrical intersection  matrix  of 
the exceptional set of $X$. Assume that  $X$ has a polygon singularity, the graph of the exceptional set is a star shaped
 with root $E_n$,  and $a_{ i,n}=1, a_{ i,i}=-2,$ for $i=1,...,n-1$. Then the matrix $A$ 
is negative definite if and only if $\displaystyle -a_{ n,n}> {{n-1}\over{2}}$ 
and $(NN)$ is true if and only if $\displaystyle -a_{ n,n}> {{n}\over{2}}$. 
So if $n$ is odd $(NN)$ is always true, but if $n$ is even it remains open the case 
$\displaystyle -a_{ n,n}= {{n}\over{2}}$. By the above theorem  only the  cases 
$(NN_{(n,i)}) $ for $i=1,...,n-1$ are not true.
\end{corollary}
\begin{example} Our theorem cannot be applied 
to the  following graph of a sandwich singularity where $a\geq 3$. In fact it follows from the theorem 
that only $(NN)_{(1,3)}$ is not true. Note that Nash problem's on arcs for (rational) sandwich 
singularities has a positive solution by the work 
of Monique Lejeune-Jalabert and Ana Reguera \cite{L-R}.
\vskip 0.5cm
 $\displaystyle A=\pmatrix{-a&0&1&0\cr 0&-2&1&0\cr 1&1&-2&1\cr 0&0&1&-2\cr}$ 
\setlength{\unitlength}{4144sp}%
\begingroup\makeatletter\ifx\SetFigFont\undefined%
\gdef\SetFigFont#1#2#3#4#5{%
  \reset@font\fontsize{#1}{#2pt}%
  \fontfamily{#3}\fontseries{#4}\fontshape{#5}%
  \selectfont}%
\fi\endgroup%
\begin{picture}(4617,0)(100,-800)
{\thinlines}%

{ \put(3286,-1140){\circle*{134}}
}%
{ \put(1846,-1140){\circle*{134}}
}%
{ \put(2566,-420){\circle*{134}}
}%
{\put(2566,-1140){\circle*{134}}}

{ \put(2566,-1140){\line( 1, 0){720}}
\put(2566,-1140){\line( -1, 0){720}}
\put(2566,-1140){\line( 0, 1){720}}
}%

\put(2450,-1330){\makebox(0,0)[lb]{\smash{{\SetFigFont{10}{16.8}{\rmdefault}{\mddefault}{\updefault}{$-2$}%
}}}}
\put(3150,-1330){\makebox(0,0)[lb]{\smash{{\SetFigFont{10}{16.8}{\rmdefault}{\mddefault}{\updefault}{$ -2$}%
}}}}
\put(1700,-1330){\makebox(0,0)[lb]{\smash{{\SetFigFont{10}{16.8}{\rmdefault}{\mddefault}{\updefault}{$ -a$}%
}}}}
\put(2450,-330){\makebox(0,0)[lb]{\smash{{\SetFigFont{10}{16.8}{\rmdefault}{\mddefault}{\updefault}{$ -2$}%
}}}}

\end{picture}%

\end{example}
\vskip 0.5cm
\begin{theorem} Let $\pi : \tilde 
X\longrightarrow X$ be the minimal resolution of singularities and let 
$A=(a_{ i,j} )$ be the $n\times n$ symmetrical intersection  matrix  of 
the exceptional set of $X$. Assume that the singularity is like a star, i.e. 
 the graph of the exceptional set is a star
 with root $E_0$ having $s\geq 3$ branches.

\setlength{\unitlength}{4144sp}%
\begingroup\makeatletter\ifx\SetFigFont\undefined%
\gdef\SetFigFont#1#2#3#4#5{%
  \reset@font\fontsize{#1}{#2pt}%
  \fontfamily{#3}\fontseries{#4}\fontshape{#5}%
  \selectfont}%
\fi\endgroup%
\begin{picture}(4617,1905)(1561,-2071)
{\thinlines
\put(1981,-781){\circle*{134}}
}%
{\put(2791,-1906){\circle*{134}}
}%
{\put(2566,-1141){\circle*{134}}
}%
{ \put(2341,-511){\circle*{134}}
}%
{ \put(1981,-1591){\circle*{134}}
}%
{ \put(2344,-1884){\circle*{134}}
}%
{ \put(3198,-1589){\circle*{134}}
}%
{
}%
{\put(3151,-781){\circle*{134}}
}%
{ \put(3286,-1140){\circle*{134}}
}%
{ \put(4029,-1140){\circle*{134}}
}%
{ \put(5424,-1140){\circle*{134}}
}%
{ \put(4724,-1140){\circle*{134}}
}%
{ \put(3871,-1996){\circle*{134}}
}%
{ \put(3759,-421){\circle*{134}}
}%

{ \put(2566,-1140){\line( 5, 3){588.971}}
}%
{ \put(2566,-1140){\line( 3,-2){675}}
}%
{ \put(2566,-1140){\line( 1,-3){238.500}}
}%
{ \put(2566,-1140){\line(-1,-3){238.500}}
}%
{ \put(2566,-1140){\line(-4,-3){590.400}}
}%
{ \put(2566,-1140){\line( 1, 0){720}}
}%
{ \put(2566,-1140){\line(-1, 3){225}}
}%
{ \put(1976,-808){\line( 5,-3){536.029}}
}%
{ \put(3331,-1140){\line( 1, 0){720}}
}%
{ \put(4051,-1141){\line( 1, 0){720}}
}%
{ \put(4726,-1141){\line( 1, 0){720}}
}%
{ \put(3234,-1601){\line( 3,-2){675}}
}%
{ \put(3207,-754){\line( 5, 3){588.971}}
}%

\put(3106,-1366){\makebox(0,0)[lb]{\smash{{\SetFigFont{10}{16.8}{\rmdefault}{\mddefault}{\updefault}{$ a(i)_{11}$}%
}}}}
\put(3841,-1366){\makebox(0,0)[lb]{\smash{{\SetFigFont{10}{16.8}{\rmdefault}{\mddefault}{\updefault}{$ a(i)_{22}$}%
}}}}
\put(4551,-1366){\makebox(0,0)[lb]{\smash{{\SetFigFont{10}{16.8}{\rmdefault}{\mddefault}{\updefault}{$ a(i)_{33}$}%
}}}}
\put(5221,-1366){\makebox(0,0)[lb]{\smash{{\SetFigFont{10}{16.8}{\rmdefault}{\mddefault}{\updefault}{$ a(i)_{44}$}%
}}}}
\put(4910,-1006){\makebox(0,0)[lb]{\smash{{\SetFigFont{10}{16.8}{\rmdefault}{\mddefault}{\updefault}{$ a(i)_{34}$}%
}}}}
\put(4170,-1006){\makebox(0,0)[lb]{\smash{{\SetFigFont{10}{16.8}{\rmdefault}{\mddefault}{\updefault}{$ a(i)_{23}$}%
}}}}
\put(3470,-1006){\makebox(0,0)[lb]{\smash{{\SetFigFont{10}{16.8}{\rmdefault}{\mddefault}{\updefault}{$ a(i)_{12}$}%
}}}}
\end{picture}%

To any branch of the star we associate a continuous fractions expansion:
$$q_i:= a(i)_{ 1,1} -{{a(i)^2_{ 1,2}}\over 
{\displaystyle a(i)_{ 2,2}-{{a(i)^2_{ 2,3}}\over {\displaystyle a(i)_{ 3,3}-{ {a(i)^2_{ 3,4}}\over {...}}}}}} $$
Then $(N N)$ is true if  we have the following condition:

\begin{itemize}
\item for any leaf  $E_i$ we have $C(A)_{i}<0$
\item  for any  vertex $E_i$ which is not the    root $C(A)_{i}\leq 0$
\item $  \forall i,j=1,...,s, i\not=j, \displaystyle  a_{0,0}+a_{0,i}+a_{0,j}-
\sum_{k=1, k\not=i,j}^{s}{{a_{0,k}^2}\over{q_{k}}}\leq 0 $
\end{itemize}

\end{theorem}
\demo Our first step consist to contract a whole branch $G_k$ of the star $\Gamma$. We reorder 
the irreducible components of the exceptional set by letting $E_n$ to be the leaf of the branch
$G_k$, $E_{n-1}$ be the unique vertex connected to $E_n$,
 $E_{n-2}$ be the unique vertex connected to $E_{n-1}$ but distinct from $E_n$, and so on until we arrive to the root named always by $E_0$,
 the order in the other branches are arbitrary. We also denote $ a _{0,k}:=a(k)_{0,1}$. By 
 contracting $E_n$ we will get again a like star graph and a new matrix $A^{(n-1)}=(a^{(n-1)}_{i,j})$ given by
 $$a^{(n-1)}_{i,j}= a^{(n)}_{i,j} - 
{{a^{(n)}_{n,i}a^{(n)}_{n,j}}\over{a^{(n)}_{n,n}} },$$ 
 regarding that our graph is a star we get:
 $$a^{(n-1)}_{i,j}= a_{i,j} {\rm \ if\ \   } (i,j)\not= (n-1,n-1){\rm \ and\ \   }$$ 
 $$a^{(n-1)}_{n-1,n-1}= a_{n-1,n-1} - 
{{{a}^2_{n,n-1}}\over{a_{n,n}}} ,$$ 
Proceeding in this way we can contract all the vertices of $G_k$, then we will get 
again a like star graph and a new matrix $A'=(a'_{i,j})$ given by
 $$a'_{i,j}= a_{i,j} {\rm \ if\ \   } (i,j)\not= (0,0){\rm \ and\ \   }$$ 
 $$a'_{0,0}= a_{0,0} - 
{{{a}^2_{0,k}}\over{q_{k}}} .$$ 
Now we are ready to  prove the claim: we need to compare any two elements
 in the graph $\Gamma$, 
these elements are in at most 
two branches $G_\alpha ,G_\beta $ of the star, we contract   $s-2$ branches (indexed by
 a set $I$) of the star
 distinct  
from $G_\alpha ,G_\beta $ 
 and we get a new graph of type $A_n$ and a  new matrix $A''=(a''_{i,j})$ given by
 $$a''_{i,j}= a_{i,j} {\rm \ if\ \   } (i,j)\not= (0,0){\rm \ and\ \   }$$ 
 $$a''_{0,0}= a_{0,0} - \sum_{k\in I}{{a_{0,k}^2}\over{q_{k}}},$$ 
the hypothesis of the theorem imply that this special tree satisfies the hypothesis of Theorem 4
 and we are done.
\section{Examples}

We discuss some examples, some of them are obtained by direct application 
of the results above. Numerical examples were computed with my software. 
Following the ideas developed in this paper
I have written a program that given in entry  the intersection matrix $A$ of
 the exceptional set
 in the minimal resolution, compute all the matrices
 $(A^\sigma)^{(l)}$ and 
check if the numerical Nash condition $(NN_{(i,j)})$  is true or not, 
the output is  a 
 $n\times n$ square matrix $N$, 
such that : $$n_{ij }=\cases{1 & if $(NN_{(i,j)})$  is true \cr
0& if $(NN_{(i,j)})$  is false \cr
\sqb & if $i=j$\cr
}.$$
 
\begin{example} Let $\pi : \tilde 
X\longrightarrow X$ be the minimal resolution of singularities and let 
$A=(a_{ i,j} )$ be the $n\times n$ symmetrical intersection  matrix  of 
the exceptional set of $X$, let $d=\max_{i\not= j}¥ \{ a_{ i,j} \}$ 
and $m=\min_{i}¥ \{ -a_{ i,i} \}$.  If $m> (n-1)d$ then the Nash map 
${ \cal N }$ is bijective.

\end{example}
It is an immediate consequence of Lemma 2. 
Remark that if   $m,d$ are  two strictly positive integers such that $m> 
(n-1)d$.  The quadratic form associated to the matrix $M$ such that 
$m_{i,i}=-m$ for any $i$, and $m_{i,j}=d$ for any $ i\not= j $ is 
negative definite.

\begin{example} If $n=3$ Nash's problem has a positive answer if 
 for any distinct numbers $i,j,k$,  we have $a_{k,i}a_{k,j} -a_{i,j}a_{k,k}< \min \{-a_{k,i}^2 +a_{i,i}a_{k,k},-a_{k,j}^2 +a_{j,j}a_{k,k} \} $. For example if $a_{i,j}\in\{0,1\}$ for any $ i\not= j$ and $-a_{1,1}\geq 2, -a_{2,2}\geq 3, -a_{3,3}\geq 3$ then the Nash's problem has a positive answer.
\end{example}
\begin{example} (NN) is  true for Rational double points $A_n$ but no true 
 for $D_n$ neither $E_6,E_7,E_8$. Remark that recently C. Plenat has proved that the Nash map 
${ \cal N }$ is bijective for the singularities $D_n.$
\begin{itemize}
\item By Proposition  2. it is enough to consider  the singularity $ D_4$, in this case we have 
 \vskip 0.5cm
 $\displaystyle A =\pmatrix{-2&1&1&1\cr 1&-2&0&0\cr 1&0&-2& 0\cr 1&0&0&-2\cr}
N=\pmatrix{\sqb&0&0&0\cr 1&\sqb&1&1\cr 1&1&\sqb&1\cr 1&1&1&\sqb\cr}$
\setlength{\unitlength}{4144sp}%
\begingroup\makeatletter\ifx\SetFigFont\undefined%
\gdef\SetFigFont#1#2#3#4#5{%
  \reset@font\fontsize{#1}{#2pt}%
  \fontfamily{#3}\fontseries{#4}\fontshape{#5}%
  \selectfont}%
\fi\endgroup%
\begin{picture}(0,0)(1500,-800)
{\thinlines}%

{ \put(3286,-1140){\circle*{134}}
}%
{ \put(1846,-1140){\circle*{134}}
}%
{ \put(2566,-420){\circle*{134}}
}%
{\put(2566,-1140){\circle*{134}}}

{ \put(2566,-1140){\line( 1, 0){720}}
\put(2566,-1140){\line( -1, 0){720}}
\put(2566,-1140){\line( 0, 1){720}}
}%

\put(2450,-1330){\makebox(0,0)[lb]{\smash{{\SetFigFont{10}{16.8}{\rmdefault}{\mddefault}{\updefault}{$-2$}%
}}}}
\put(3150,-1330){\makebox(0,0)[lb]{\smash{{\SetFigFont{10}{16.8}{\rmdefault}{\mddefault}{\updefault}{$ -2$}%
}}}}
\put(1700,-1330){\makebox(0,0)[lb]{\smash{{\SetFigFont{10}{16.8}{\rmdefault}{\mddefault}{\updefault}{$ -2$}%
}}}}
\put(2450,-330){\makebox(0,0)[lb]{\smash{{\SetFigFont{10}{16.8}{\rmdefault}{\mddefault}{\updefault}{$ -2$}%
}}}}
\end{picture}%
\vskip 0.5cm
\item Consider the singularity $ E_6$, in this case we have 
$$A=\pmatrix{-2&1&0&0&0&0\cr 1&-2&1&0&0&0\cr 0&1&-2&1&0&1\cr 0&0&1&-2&1&0\cr
0&0&0&1&-2&0\cr0&0&1&0&0&-2\cr}
N=\pmatrix{\sqb&1&1&1&1&1\cr 0&\sqb&1&1&0&0\cr 0&0&\sqb&0&0&0\cr 0&1&1&\sqb&0&0\cr 1&1&1&1&\sqb&1\cr 1&1&1&1&1&\sqb\cr}
$$
\item Consider the singularity $ E_7$, in this case we have 

$$A=\pmatrix{-2&1&0&0&0&0&0\cr 1&-2&1&0&0&0&0\cr 0&1&-2&1&0&0&1\cr
 0&0&1&-2&1&0&0\cr 0&0&0&1&-2&1&0\cr 0&0&0&0&1&-2&0\cr 0&0&1&0&0&0&-2\cr}
N=\pmatrix{\sqb&1&1&1&1&1&1\cr 0&\sqb&1&1&0&0&0\cr 0&0&\sqb&0&0&0&0\cr 
0&0&1&\sqb&0&0&0\cr 0&1&1&1&\sqb&0&1\cr 1&1&1&1&1&\sqb&1\cr
 0&1&1&1&1&0&\sqb\cr}$$
\item Consider the singularity $ E_8$, in this case we have 

$$A=\pmatrix{-2&1&0&0&0&0&0&0\cr 1&-2&1&0&0&0&0&0\cr 0&1&-2&1&0&0&0&1\cr
 0&0&1&-2&1&0&0&0\cr 0&0&0&1&-2&1&0&0\cr 0&0&0&0&1&-2&1&0\cr 0&0&0&0&0&1&-2&0\cr 0&0&1&0&0&0&0&-2\cr}
N=\pmatrix{\sqb&1&1&1&1&1&0&1\cr 0&\sqb&1&1&0&0&0&0\cr 0&0&\sqb&0&0&0&0&0\cr 0&0&1&\sqb&0&0&0&0\cr 
0&1&1&1&\sqb&0&0&0\cr 0&1&1&1&1&\sqb&0&1\cr 1&1&1&1&1&1&\sqb&1
\cr0&1&1&1&1&0&0&\sqb\cr }$$

\end{itemize}
\end{example}

\begin{example} The following two graphs are like star shaped, and condition $(NN)$ is not true. 
 $$A=\pmatrix{-2&1&0&0&0&0&0&0&0&0\cr 1&-2&1&0&0&0&0&0&0&0\cr 0&1&-3&1&0&1&0&0&1&0\cr
0&0&1&-2&1&0&0&0&0&0\cr 0&0&0&1&-2&0&0&0&0&0\cr 0&0&1&0&0&-2&1&0&0&0\cr
 0&0&0&0&0&1&-2&1&0&0\cr 0&0&0&0&0&0&1&-2&0&0\cr 0&0&1&0&0&0&0&0&-2&1\cr
 0&0&0&0&0&0&0&0&1&-2\cr
 }
N=\pmatrix{\sqb&1&1&1&1&1&1&1&1&1\cr 0&\sqb&1&1&1&1&1&1&1&1
\cr 0&0&\sqb&0&0&0&1&0&0&0\cr
 1&1&1&\sqb&0&1&1&1&1&1\cr 1&1&1&1&\sqb&1&1&1&1&1\cr 1&1&1&1&1&\sqb&0&0&1&1\cr
 1&1&1&1&1&1&\sqb&0&1&1\cr 1&1&1&1&1&1&1&\sqb&1&1\cr 1&1&1&1&1&1&1&1&\sqb&0\cr 1&1&1&1&1&1&1&1&1&\sqb\cr }$$
  $$A=\pmatrix{-2 & 1 & 0 & 0 & 0 & 0 & 0 & 0 & 0\cr 1 & -3 & 1 & 0 & 0 & 1 & 0 & 1 & 0\cr
 0 & 1 & -2 & 1 & 0 & 0 & 0 & 0 & 0\cr 0 & 0 & 1 & -2 & 1 & 0 & 0 & 0 & 0\cr
 0 & 0 & 0 & 1 & -2 & 0 & 0 & 0 & 0\cr 0 & 1 & 0 & 0 & 0 & -2 & 1 & 0 & 0\cr
 0 & 0 & 0 & 0 & 0 & 1 & -2 & 0 & 0\cr  0 & 1 & 0 & 0 & 0 & 0 & 0 & -2 & 1\cr
 0 & 0 & 0 & 0 & 0 & 0 & 0 & 1 & -2\cr}
N=\pmatrix{\sqb & 1 & 1 & 1 & 1 & 1 & 1 & 1 & 1\cr 0 & \sqb & 1 & 1 & 1 & 1 & 1 & 1 & 1\cr
 1 & 1 & \sqb & 1 & 1 & 1 & 1 & 1 & 1\cr 1 & 1 & 1 & \sqb & 1 & 1 & 1 & 1 & 1\cr
 1 & 1 & 1 & 1 & \sqb & 1 & 1 & 1 & 1\cr 1 & 1 & 1 & 1 & 1 & \sqb & 1 & 1 & 1\cr
 1 & 1 & 1 & 1 & 1 & 1 & \sqb & 1 & 1\cr  1 & 1 & 1 & 1 & 1 & 1 & 1 & \sqb & 1\cr
 1 & 1 & 1 & 1 & 1 & 1 & 1 & 1 & \sqb\cr}$$
\end{example}
\begin{example} The following  graphs are like star shaped, and condition $(NN)$ is  true.
 $$A=\pmatrix{-2 & 1 & 0 & 0 & 0 & 0 & 0 & 0\cr 1 & -3 & 1 & 0 & 0 & 1 & 0 & 1\cr 
0 & 1 & -2 & 1 & 0 & 0 & 0 & 0\cr 0 & 0 & 1 & -2 & 1 & 0 & 0 & 0\cr
 0 & 0 & 0 & 1 & -2 & 0 & 0 & 0\cr 0 & 1 & 0 & 0 & 0 & -2 & 1 & 0\cr 
 0 & 0 & 0 & 0 & 0 & 1 & -2 & 0\cr
0 & 1 & 0 & 0 & 0 & 0 & 0 & -2\cr}
N=\pmatrix{\sqb & 1 & 1 & 1 & 1 & 1 & 1 & 1\cr 1 & \sqb & 1 & 1 & 1 & 1 & 1 & 1\cr
 1 & 1 & \sqb & 1 & 1 & 1 & 1 & 1\cr 1 & 1 & 1 & \sqb & 1 & 1 & 1 & 1\cr 
1 & 1 & 1 & 1 & \sqb & 1 & 1 & 1\cr 
1 & 1 & 1 & 1 & 1 & \sqb & 1 & 1\cr 1 & 1 & 1 & 1 & 1 & 1 & \sqb & 1\cr
1 & 1 & 1 & 1 & 1 & 1 & 1 & \sqb\cr}$$

\end{example}
\begin{example}In this example the graph of the singularity is a tree, $(NN)$ is true but we can't apply  Theorem 4.

$$A=\pmatrix{-2&1&0&0&0&0&0&0\cr 1&-3&1&0&1&0&0&0\cr 0&1&-2&0&0&0&0&0\cr 
0&0&0&-2&1&0&0&0\cr 0&1&0&1&-3&1&0&1\cr 0&0&0&0&1&-2&1&0\cr 0&0&0&0&0&1&-2&0\cr
 0&0&0&0&1&0&0&-2\cr }
N=\pmatrix{\sqb & 1 & 1 & 1 & 1 & 1 & 1 & 1\cr 1 & \sqb & 1 & 1 & 1 & 1 & 1 & 1\cr
 1 & 1 & \sqb & 1 & 1 & 1 & 1 & 1\cr 1 & 1 & 1 & \sqb & 1 & 1 & 1 & 1\cr
 1 & 1 & 1 & 1 & \sqb & 1 & 1 & 1\cr 1 & 1 & 1 & 1 & 1 & \sqb & 1 & 1\cr
 1 & 1 & 1 & 1 & 1 & 1 & \sqb & 1\cr 1 & 1 & 1 & 1 & 1 & 1 & 1& \sqb\cr }$$
\end{example}

{\bf Acknowledgment}. The author thanks Camille Plenat and the referee for helpful comments.

\end{document}